\documentclass[letterpaper, 10 pt, conference]{ieeeconf}  

\IEEEoverridecommandlockouts                              

\title{\LARGE \bf Two-Level Decentralized-Centralized Control of Distributed Energy Resources in Grid‐Interactive Efficient Buildings 
}

\author{Xiang Huo$^{*}$, Jin Dong$^{\dagger}$, Borui Cui$^{\dagger}$, Boming Liu$^{\dagger}$, Jianming Lian$^{\dagger}$, and Mingxi Liu$^{*}$
\thanks{This manuscript has been authored by UT-Battelle, LLC, under contract DE-AC05-00OR22725 with the US Department of Energy (DOE). The US government retains and the publisher, by accepting the article for publication, acknowledges that the US government retains a nonexclusive, paid-up, irrevocable, worldwide license to publish or reproduce the published form of this manuscript, or allow others to do so, for US government purposes. DOE will provide public access to these results of federally sponsored research in accordance with the DOE Public Access Plan \href{http://energy.gov/downloads/doe-public-access-plan}.}
\thanks{$^{*}$Xiang Huo and Mingxi Liu are with the Department of Electrical and Computer Engineering,
University of Utah,
Salt Lake City, UT, USA
        {\tt\small \{xiang.huo,mingxi.liu\}@utah.edu}}%
\thanks{$^{\dagger}$Jin Dong, Borui Cui,
Boming Liu, and Jianming Lian are with the Electrification and Energy Infrastructures Division, Oak Ridge National Laboratory, Oak Ridge, TN, USA
        {\tt\small \{dongj,cuib,liub, lianj\}@ornl.gov}}%
\thanks{\noindent This work has been supported in part by NSF Award: ECCS-2145408.}
}

\usepackage{cite}
\usepackage{amsmath,amssymb,amsfonts}
\usepackage{graphicx}
\usepackage{textcomp}

\usepackage{bm}

\usepackage[linesnumbered]{algorithm2e}

\usepackage{amssymb}

\usepackage{subcaption}

\newcounter{defcounter}
\newenvironment{pequation}{%
\addtocounter{equation}{-1}
\refstepcounter{defcounter}
\renewcommand\theequation{{\bf{P\thedefcounter}}}
\begin{equation}}
{\end{equation}}

\usepackage[hang,flushmargin]{footmisc}
\usepackage{lipsum}
\makeatletter
\newcommand{\algorithmfootnote}[2][\footnotesize]{%
  \let\old@algocf@finish\@algocf@finish
  \def\@algocf@finish{\old@algocf@finish
    \leavevmode\rlap{\begin{minipage}{\linewidth}
    #1#2
    \end{minipage}}%
  }%
}
\makeatother

\usepackage{soul}
\usepackage{xargs}

\usepackage[pdftex,dvipsnames]{xcolor}  
\usepackage[colorinlistoftodos,prependcaption,textsize=scriptsize, textwidth=8mm]{todonotes}
\newcommandx{\revfirst}[2][1=]{\todo[linecolor=blue,backgroundcolor=blue!25,bordercolor=blue,#1]{#2}}
\newcommandx{\revsecond}[2][1=]{\todo[linecolor=red,backgroundcolor=red!25,bordercolor=red,#1]{#2}}
\newcommandx{\revthird}[2][1=]{\todo[linecolor=ForestGreen,backgroundcolor=ForestGreen!25,bordercolor=ForestGreen,#1]{#2}}
\newcommandx{\revfourth}[2][1=]{\todo[linecolor=OliveGreen,backgroundcolor=OliveGreen!25,bordercolor=OliveGreen,#1]{#2}}
\newcommandx{\revfifth}[2][1=]{\todo[linecolor=orange,backgroundcolor=orange!25,bordercolor=orange,#1]{#2}}

\newcommandx{\thiswillnotshow}[2][1=]{\todo[disable,#1]{#2}}

\begin{document}

\maketitle
\thispagestyle{empty}
\pagestyle{empty}


\begin{abstract}

The flexible, efficient, and reliable operation of  grid-interactive efficient buildings (GEBs) is increasingly impacted by the growing penetration of distributed energy resources (DERs). Besides, the optimization and control of DERs, buildings, and distribution networks are further complicated by their interconnections. In this paper, we exploit load-side flexibility and clean energy resources to develop a novel two-level hybrid decentralized-centralized (HDC) algorithm to control DER-connected GEBs. The proposed HDC 1) achieves scalability \textit{w.r.t.} to a large number of grid-connected buildings and devices,  2) incorporates a two-level design where aggregators control buildings centrally and the system operator coordinates the distribution network in a decentralized fashion, and 3) improves the computing efficiency and enhances communicating compatibility with heterogeneous temporal scales. Simulations are conducted based on the prototype of a campus building at the {Oak Ridge National Laboratory} to show the efficiency and efficacy of the proposed approach.

\end{abstract}

\section{Introduction}

Grid-interactive efficient buildings (GEBs) aim at revolutionizing traditional buildings into clean and flexible energy assets by integrating the ever-growing distributed energy resources (DERs) and flexible electric loads \cite{grid_building}.  
The rapid employment of DERs, such as solar photovoltaics (PVs), electric vehicles (EVs), and energy storage systems (ESSs), can significantly accelerate building electrification, improve grid resilience, decrease carbon emissions, and reduce infrastructure cost \cite{zhang2018distributed,akorede2010distributed}. Besides, the substantial flexible electric loads in buildings such as heating, ventilation, and air conditioning (HVAC) can be controlled to serve customers' needs and maximize the building's energy efficiency.  
Therefore, we expect to leverage the potential of DER-connected GEBs with optimized solutions to provide both grid-level and customer-side services. 

Controlling grid-edge resources (e.g., DERs and HVACs) in a centralized fashion is easy to implement \cite{majumdar2017centralized,ciocia2018voltage,zhang2019distributionally}, but they inevitably suffer from poor reliability and scalability. To overcome the drawbacks, distributed  methods can achieve higher self-determination by assigning the tasks to the agents. 
Fan \emph{et al.} in \cite{fan2021distributed} 
proposed a distributed discrete-time control scheme to achieve the optimal coordination of conventional and renewable generators. In \cite{zhang2021distributed}, an asynchronous distributed leader-follower control method was proposed to achieve conservation voltage reduction by optimally scheduling smart inverters of DERs. 
Wang \emph{et al.} in  \cite{wang2019distributed} proposed a hierarchical distributed scheme that utilizes thermostatically controlled loads to provide ancillary services. Though distributed methods bring scalability, the frequent peer-to-peer 
communications can
limit their applications. 

In contrast, decentralized strategies require no communication between agents or subsystems, which can significantly ease the communication burdens. In \cite{lin2017decentralized}, a decentralized disturbance-feedback controller was designed to coordinate and control PV inverters and ESSs. 
Kou \emph{et al.} in \cite{xiao2020ADMM,kou2020comprehensive} developed alternating direction method of multipliers based algorithms to coordinate residential demand-side resources and ensure efficient operation of distribution networks. In \cite{doghmane2020design}, a high dimensional decentralized controller was proposed to improve the robustness and optimality in the application of a smart building system. 

Despite the advantages of decentralized methods, few state-of-the-art research has addressed the decentralized control of coupled objectives and constraints in terms of buildings, DERs, and distribution networks. Besides, the computing burden imposed on individual agents and the synchronization requirement on iterative updates further restrict the practicality of decentralized methods. To resolve those issues, we expect to design a hybrid decentralized-centralized (HDC) algorithm that allows  a building aggregator to control its building and the system operator (SO) of the distribution network to coordinate asynchronous updates in a decentralized way with heterogeneous temporal scales. 

The contributions of this paper are four-fold: 1) A novel two-level HDC strategy is proposed where aggregators  centrally control DERs and HVACs at the building level and the SO coordinates at the distribution network level in a decentralized fashion; 2) The proposed method achieves scalability \emph{w.r.t.} the number of DERs as well as GEBs in a distribution network; 3) Asynchronous communication with heterogeneous temporal scales is investigated to improve the system compatibility and computing efficiency; 4) We benchmark the problem formulation and algorithm design considering both local and global objectives and constraints, and verify the efficiency and efficacy of the proposed method via a real-world setup.

\section{Mathematical Formulation}

\subsection{System Modeling} 

\subsubsection{Building envelop thermal network model} 
\begin{figure}[!htb]
    \centering
\includegraphics[width=0.38\textwidth, trim={0cm 0cm 0.0cm 0cm},clip]{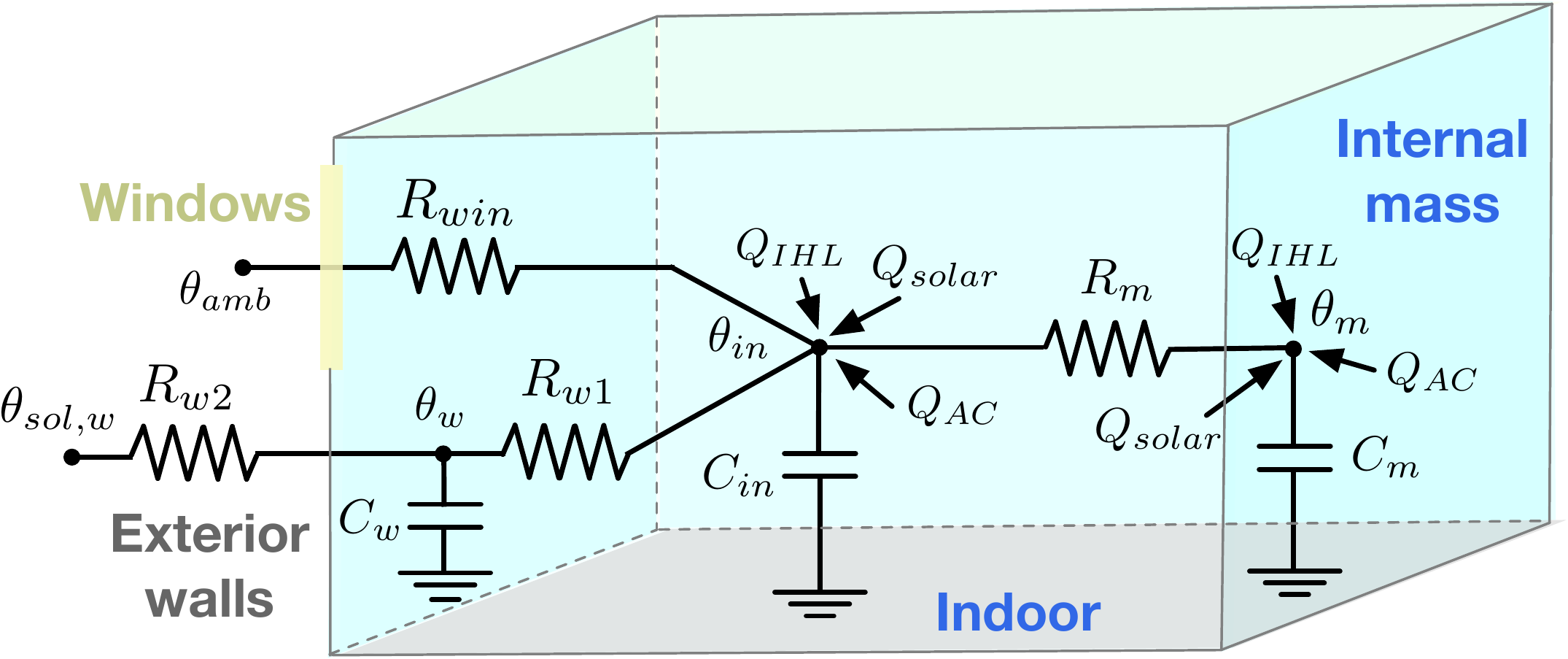}
\caption{The thermal network model of an office room (room 102 inside building 3147 at ORNL).\label{building_thermal_network}}
\end{figure}

Consider a single room in building 3147 at the {Oak Ridge National Laboratory} (ORNL). The heat transfer of this room can be described by the resistance-capacitance (RC) thermal network model \cite{cui2019load} as
\begin{align}
C_{w} \frac{d\theta_{w}(t)}{dt} &=\frac{\theta_{sol,w}(t)-\theta_{w}(t)}{R_{w2}}-\frac{\theta_{w}(t)-\theta_{in}(t)}{R_{w1}}   
\label{1}\\
C_{in} \frac{d\theta_{in}(t)}{d t}&=\frac{\theta_{w}(t)-\theta_{i n}(t)}{R_{w1}}+\frac{\theta_{am b}(t)-\theta_{i n}(t)}{R_{win}}\nonumber\\
&\quad +\frac{\theta_{m}(t)-\theta_{in}(t)}{R_{m}}+Sp_{1}Q_{AC}u(t)\nonumber\\
&\quad +Sp_{2}Q_{IHL} +Sp_{3}Q_{solar}
\label{2}\\
    C_{m} \frac{d\theta_{m}(t)}{d t}&=-\frac{\theta_{m}(t)-\theta_{in}(t)}{R_{m}}+\left(1-Sp_{1}\right)Q_{AC}u(t)\nonumber\\
    &\quad+\left(1-Sp_{2}\right)Q_{IHL}+\left(1-Sp_{3}\right)Q_{solar}
    \label{3}
\end{align} where $C_{w}$, $C_{in}$, and $C_{m}$ denote the thermal capacitances of the exterior wall, indoor air, and internal mass, respectively. $R_{w1}$ ($R_{w2}$), $R_{win}$, and $R_m$ denote the thermal resistances of exterior walls, window, and internal mass, respectively. $\theta_{w}(t)$, $\theta_{in}(t)$, $\theta_{m}(t)$, $\theta_{amb}(t)$, and $\theta_{sol,w}(t)$ are the envelope temperature, indoor temperature, indoor thermal mass, outdoor dry bulb temperature, and solar air temperature on the external surface of building envelope at time $t$, respectively. $Sp_{1}$, $Sp_{2}$, and $Sp_{3}$ denote the convection fractions. Let $u(t) \in \{0,1\}$ denote the \textsc{On-Off} status and $Q_{AC}$ denote cooling capacity of the HVAC, respectively, and $Q_{IHL}$, and $Q_{solar}$ denote indoor heat load and the solar radiation through windows, respectively. The 3R3C thermal network model of the office room is shown in Fig. \ref{building_thermal_network}.

In the 3R3C model, the discrete-time dynamics of $\theta_{w}(t)$, $\theta_{in}(t)$, and $\theta_{m}(t)$ with the sampling time $\Delta T$ can be represented as 
\begin{equation}
    \theta(t+1) = \Delta T \frac{d\theta(t)}{dt} + \theta(t).
    \label{4}
\end{equation}

To guarantee the indoor temperature stay within the residents' comfort zone, $\bm{\theta}_{in}$ should satisfy 
\begin{equation}
   \hat{\bm{\theta}}_{l} \leq   \bm{\theta}_{in} \leq \hat{\bm{\theta}}_{u}
   \label{5}
\end{equation}
where $\bm{\theta}_{in} = [\theta_{in}(1), \ldots , \theta_{in}(T)]^\mathsf{T}$, $T$ denotes the final time interval, $\hat{\bm{\theta}}_{l}$ and $\hat{\bm{\theta}}_{u}$ denote the lower and upper temperature bounds, respectively. 


\noindent \textbf{Remark 1:} Eqs. \eqref{1}-\eqref{5} describe the heat transfer model and temperature constraint of a single room. Similarly, by constructing the 3R3C models of all rooms, the indoor temperatures of individual rooms can be precisely controlled. However, the inclusion of all room models can increase the modeling complexity and computing overhead. To improve the computing efficiency, this paper adopts the aggregator model which aggregates heterogeneous RC parameters based on the individual HVACs \cite{wang2021tri} without sacrificing accuracy. More details can be found in Section \ref{building_centralize}. \hfill $\square$



\subsubsection{Electric vehicle}

Let $\bm{p}_{v\epsilon} \in \mathbb{R}^{T}$ denote the charging profile of the $\epsilon$th EV, and it is constrained by 
\begin{equation}
    \bm{0} \leq \bm{p}_{v\epsilon} \leq \bm{r}_v^u \label{1s}
\end{equation}
where $\bm{r}_v^u$ denotes the maximum charging power. To guarantee all EVs can be charged to the desired energy level, the total charging loads of the $\epsilon$th EV should satisfy
\begin{equation}
    \bm{G}\bm{p}_{v\epsilon} = d_\epsilon \label{2s}
\end{equation}
where $\bm{G} \triangleq [\Delta T,\ldots,\Delta T] \in \mathbb{R}^{1\times T}$ denotes the aggregation vector and $d_\epsilon$ denotes the charging demand of the $\epsilon$th EV.

\subsubsection{Solar photovoltaic} 
During $T$ time intervals of a day, the active power injections from the $k$th PV are limited by
\begin{equation}
\bm{0} \leq \bm{p}_{sk} \leq \bm{\hat{p}}^{u}_{s}
\label{pv_limit}
\end{equation}
where $\bm{p}_{sk} \in \mathbb{R}^{T}$ denotes the active power injections and $\bm{\hat{p}}^{u}_{s}$ denotes the maximum available active power from the $k$th PV inverter. $\bm{\hat{p}}_{s}^{u}$ is assumed to be known by forecast.

\subsubsection{Energy storage system} The discharging/charging power of the $\sigma$th ESS should satisfy
\begin{equation}
-\hat{\bm{p}}_{e}^{l} \leq \bm{p}_{e\sigma} \leq \hat{\bm{p}}_{e}^{u}
\label{charge_limit}
\end{equation}
where $\bm{p}_{e\sigma} \in \mathbb{R}^{T}$ denotes the discharging/charging power of the $\sigma$th ESS, and $\hat{\bm{p}}_{e}^{l}$ and $\hat{\bm{p}}_{e}^{u}$ denote the maximum discharging and charging power, respectively.
The capacity limit of the $\sigma$th ESS is constrained by 
\begin{equation}
 \hat{\bm{p}}_e^{cl} \leq \bm{A} \bm{p}_{e\sigma} \Delta T \leq
 \hat{\bm{p}}_e^{cu} \label{state_limit}
\end{equation}
where $\hat{\bm{p}}_e^{cl}$ and $\hat{\bm{p}}_e^{cu}$ denote the lower and upper capacity bounds of the ESSs, respectively, and $\bm{A}$ is a lower-triangular aggregation matrix, i.e., $\bm{A}_{\hat{\imath},\hat{\jmath}} = 1 ~ \text{if}~ \hat{\imath} \geq \hat{\jmath},  \bm{A}_{\hat{\imath},\hat{\jmath}} = 0 ~ \text{if}~ \hat{\imath} < \hat{\jmath}, \forall \hat{\imath},\hat{\jmath} = 1,\ldots,T$.


\subsubsection{Distribution network model}
\begin{figure}[!htb]
\vspace{-1mm}
    \centering
    \includegraphics[width=0.45\textwidth, trim={0.1cm 0cm 0.0cm 0cm},clip]{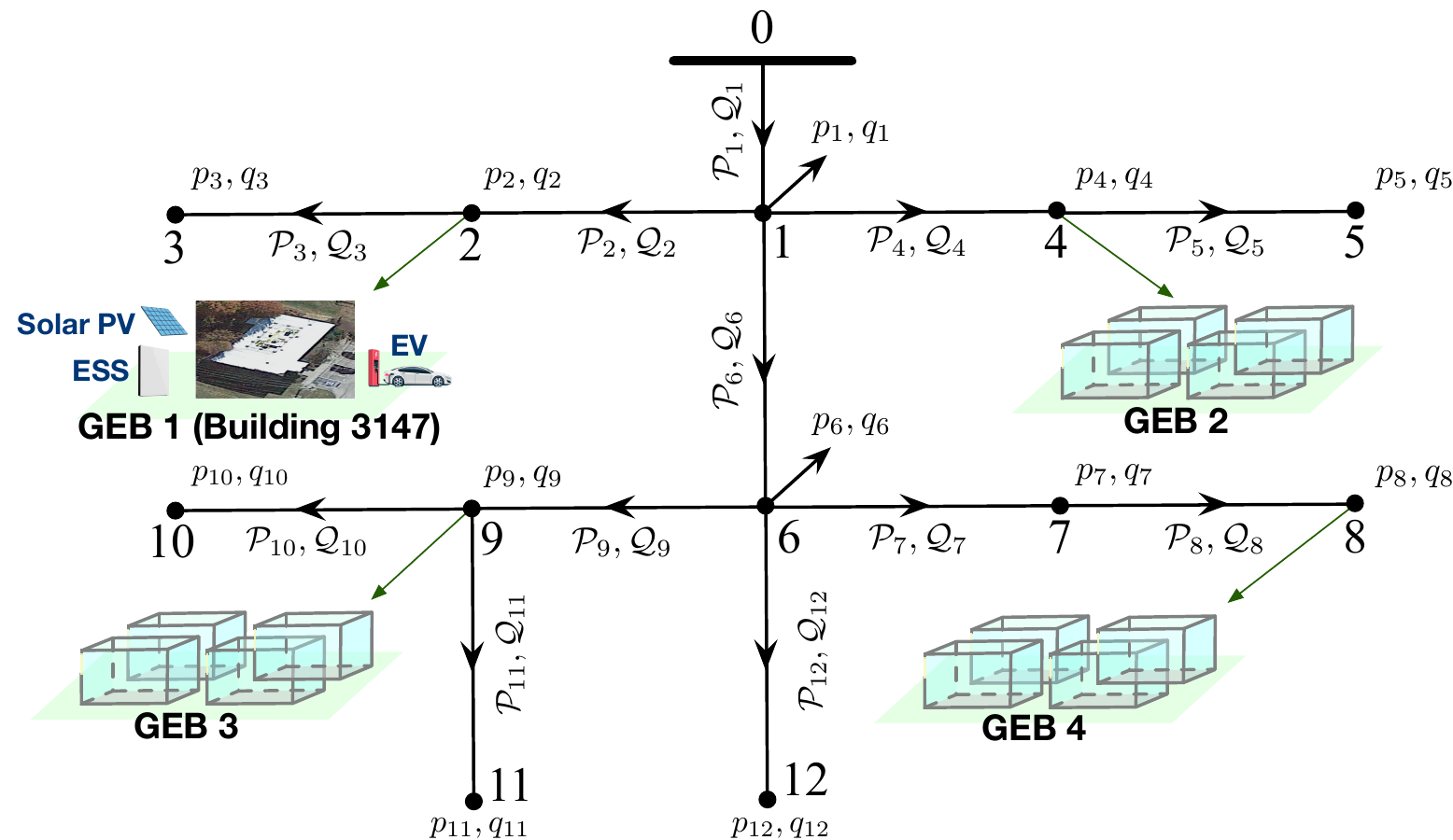}
    \caption{Integration of GEBs into the 13-node distribution network.}
\label{buildings_distribution_network}
     \vspace{-4mm}
\end{figure}

 To demonstrate the grid-level applications, we consider the control of GEBs in distribution networks. The integration of four GEBs into a 13-node distribution network is shown in Fig. \ref{buildings_distribution_network} where GEB 1 is represented by building 3147 at ORNL.  Consider a re-indexed radial distribution network and define $\mathbb{N}=\{i \mid i=1, \ldots, n\}$ as the set of downstream nodes. Let $\mathbb{L}$ denote the set of all downstream line segments and $l_{ij}$ denote the line connecting node $i$ and node $j$. Following the linear DistFlow branch equations \cite{baran1989network}, the voltage magnitude at node $i$ across $T$ time intervals is described by
\begin{equation} \label{11sss}
\boldsymbol{V}_i =\boldsymbol{V}_{0}- 2\sum_{\iota=1}^{n} \boldsymbol{R}_{i\iota} \boldsymbol{p}_\iota - 2\sum_{\iota=1}^{n} \boldsymbol{X}_{i\iota} \boldsymbol{q}_\iota
\end{equation}
where $\boldsymbol{V}_i {=}[\left|V_{i}(1)\right|^{2}, \ldots, \left|V_{i}(T)\right|^{2}]^{\mathsf{T}} \in \mathbb{R}^{T}$ and $|V_i(t)|$ is the nodal voltage magnitude, $\boldsymbol{V}_{0} =[\left|V_{0}\right|^{2}, \ldots, \left|V_{0}\right|^{2}]^{\mathsf{T}} \in \mathbb{R}^{T}$, $|V_0|$ denotes the voltage magnitude of the feeder head, $\boldsymbol{p}_i {=}[p_{i}(1), \ldots, p_{i}(T)]^{\mathsf{T}} \in \mathbb{R}^{T} $ and $\boldsymbol{q}_i {=}\left[q_{i}(1),\ldots, q_{i}(T)\right]^{\mathsf{T}} \in \mathbb{R}^{T}$ denote the active and reactive loads at node $i$, respectively, and the adjacency matrices $\bm{R}$ and $\bm{X}$ are defined by
\begin{align}
\boldsymbol{R} &\in \mathbb{R}^{n \times n},  \boldsymbol{R}_{ij}=\sum_{(\hat{\imath}, \hat{\jmath}) \in \mathbb{L}_{i} \cap \mathbb{L}_{j}} r_{\hat{\imath} \hat{\jmath}} \nonumber\\ 
\boldsymbol{X} &\in \mathbb{R}^{n \times n}, \boldsymbol{X}_{ij}=\sum_{(\hat{\imath}, \hat{\jmath}) \in \mathbb{L}_{i} \cap \mathbb{L}_{j}} x_{\hat{\imath} \hat{\jmath}}\nonumber
\end{align} 
where $r_{\hat{\imath} \hat{\jmath}}$ and $x_{\hat{\imath} \hat{\jmath}}$ denote the  resistance and reactance of line $l_{ij}$, respectively. $\mathbb{L}_{i} \subseteq \mathbb{L}$ is the set containing downstream line segments connecting the feeder head and node $i$. The nodal voltage should remain within the voltage limit
\begin{equation}
\hat{\bm{V}}_{l} \leq \boldsymbol{V}_i \leq \hat{\bm{V}}_{u}, \forall i = 1,\ldots,n \label{13}
\end{equation} 
where $\hat{\bm{V}}_{l}$ and $\hat{\bm{V}}_{u}$ denote the lower and upper voltage bounds, respectively. 



\subsection{Objective Functions} 

\subsubsection{Power loss minimization}



The global objective aims at minimizing the total active power loss in the distribution network, it can be calculated by
\begin{equation}
    f_1(\bm{p}) = \sum_{l_{ij} \in \mathbb{L}} r_{ij} \left
(\frac{\| \bm{\mathcal{P}}_{ij} \|^2_2 + \| \bm{\mathcal{Q}}_{ij} \|^2_2}{|V_0|^2}\right) \label{14}
\end{equation}
where $\bm{\mathcal{P}}_{ij} \in \mathbb{R}^{T}$ and 
$\bm{\mathcal{Q}}_{ij} \in \mathbb{R}^{T}$ denote the active and reactive power flow from node $i$ to node $j$, respectively. We assume $\bm{\mathcal{Q}}_{ij}$'s to be known constant vectors.

\subsubsection{Energy efficiency cost} For energy efficiency objectives, we only consider the ESS's degradation cost as a representative for clarity. Other cost functions can be referred to \textbf{Remark 2}. The $\sigma$th ESS's degradation cost can be calculated in terms of the smoothness of charging and discharging by
\begin{equation}
    f_2(\bm{p}_{e\sigma}) = \| \bm{B}\bm{p}_{e\sigma}  \|^2_2.
    \label{3s}
\end{equation}
where $\bm{B}$ aggregates discharging/charging differences between adjacent times, i.e., $\bm{B}_{\hat{\imath},\hat{\imath}} = 1$, $\forall \hat{\imath} = 1,\ldots,T$, $\bm{B}_{\hat{\imath},\hat{\imath}+1} = -1, \forall \hat{\imath} = 1,\ldots,T-1$, and all other elements are zeros. 


Therefore, the GEB control problem is formulated into a constrained optimization problem that aims at minimizing active power loss and improving energy efficiency as
\begin{pequation} \label{p1}
\begin{aligned}
& \underset{\bm{p},\bm{u}}{\min} & &   \delta_1 f_1(\bm{p})  {+} \sum_{\sigma=1}^{\mathcal{E}}\delta_2f_2(\bm{p}_{e\sigma})
 \\
& \: \text{s.t.} & & \eqref{1}-\eqref{13}
\end{aligned}
\end{pequation}where $\delta_1$ and $\delta_2$ are weighting constants. $\bm{p} = [\bm{p}_1, \ldots, \bm{p}_n]$, $\bm{u} = [\bm{u}_1, \ldots, \bm{u}_m]$, $m$ and $\mathcal{E}$ denote total number of GEBs and ESSs, respectively,  $\bm{u}_j$ denotes the \textsc{On-Off} decision variables of all HVACs in building $j$ accross $T$ time intervals.

\noindent \textbf{Remark 2:} Problem \eqref{p1} aims at maximizing the benefits for both customers and the SO, i.e., both global and local objectives (constraints) are considered. Depending on the practical applications, extensions of local objectives, such as the curtailment cost of the $k$th solar PV, i.e., $ \|\bm{p}_{sk} - \bm{\hat{p}}^{u}_{s}\|^2_2$, the utility cost of all the nodes, i.e., $\sum_{i=1}^{n} \bm{p}_i^{\mathsf{T}} \bm{\chi}$ where $\bm{\chi} \in \mathbb{R}^{T}$ denotes the daily real-time electricity prices, or global power flow limit constraints, can be readily incorporated. \hfill $\square$

\noindent \textbf{Remark 3:} In problem \eqref{p1}, two objectives were included in the objective function, i.e., $f_1(\bm{p})$ as global objective representing total power loss and $f_2(\bm{p}_{e\sigma})$ as local objective representing energy efficiency cost. Having $\delta_1$ and $\delta_2$ before $f_1(\bm{p})$ and $f_2(\bm{p}_{e\sigma})$ can 1) allow different weightings based on the stakeholders' desired rewards from global and local objective functions and 2) regulate different units, e.g., power loss is in kW and ESS degradation cost is in (kW)$^2$. Therefore, stakeholders can select the values of $\delta_1$ and $\delta_2$ to deal with diverse practical applications.  \hfill $\square$

\section{Algorithm Design}

In this section, we propose a scalable two-level HDC algorithm to improve the computing efficiency and compatibility of GEB control. The DERs and HVACs are centrally controlled at the building level through an aggregator while GEBs are coordinated by the SO in a decentralized way at the grid level. The SO and aggregators are assumed to be able to communicate bilaterally. 

\subsection{Building-level Centralized Control} \label{building_centralize}

{\color{blue}}
\subsubsection{Aggregator model} The \textsc{On-Off} operation of the HVACs inside a building is centrally controlled by an aggregator \cite{wang2021tri}. The aggregator integrates the parameters of the 3R3C model from all individual rooms and obtains the following first-order temperature dynamic equation as 
\begin{equation}
\tilde{\theta}_{in,j}(t+1)=a_{j} \tilde{\theta}_{in,j}(t)+b_{j} \theta_{amb}(t)+g_{j} \tilde{u}_j(t) \label{21s}
\end{equation}
where $a_{j}$, $b_{j}$, and $g_{j}$ denote the aggregated parameters, $\tilde{\theta}_{in,j}(t)$ denotes the average temperature of building $j$, $\tilde{u}_j(t)$ is a continuous variable denoting the percentage of \textsc{On}-state  HVACs in the building and is constrained by 
\begin{equation}
0 \leq \tilde{u}_j(t) \leq 1. \label{22s}
\end{equation}
The average temperature of the building should stay within the comfort zone as 
\begin{equation}
   \hat{\bm{\theta}}_{l} \leq   \tilde{\bm{\theta}}_{in,j} \leq \hat{\bm{\theta}}_{u}
   \label{17s}
\end{equation}
where $\tilde{\bm{\theta}}_{in,j} = [\tilde{\theta}_{in,j}(1), \ldots, \tilde{\theta}_{in,j}(T)]^\mathsf{T}$. The aggregated active power consumption of the HVACs in the $j$th building can be calculated by 
\begin{equation}
\bm{p}_{hj} = N_{j} P^{r} \tilde{\bm{u}}_j \label{23s}
\end{equation}
where $\bm{p}_{hj}$ denotes the aggregated active power consumption, $\tilde{\bm{u}}_j = [\tilde{\bm{u}}_j(1), \ldots, \tilde{\bm{u}}_j(T)]^\mathsf{T}$,  $N_{j}$ denotes the total number of HVACs, and $P^{r}$ denotes the HVAC's rated power. 

After obtaining the aggregated decision variable $\tilde{\bm{u}}_j$, the aggregator can first calculate the total number of \textsc{On}-state HVACs by $N_j\tilde{\bm{u}}_j$, then determine the specific \textsc{On}-state devices based on priority list strategy, e.g., based on the descending order of temperatures deviations from the temperature set point that can be calculated as $\| \tilde{\bm{\theta}}_{in,j} {-} (\hat{\bm{\theta}}_{l} {+} \hat{\bm{\theta}}_{u})/2 \|_2$. The priority list can also include the response delays or resource utilization to achieve specific building-level control goals \cite{lu2012design,wu2017hierarchical}. For the control of local DERs, the aggregator can simply adopt the updated decision variables.

\subsection{Grid-level Decentralized Coordination}
By adopting the aggregator model to centrally control each building, problem \eqref{p1} can be reformulated as  \vspace{-2mm}
\begin{pequation} \label{p2}
\begin{aligned}
& \underset{\bm{p},\tilde{\bm{u}}}{\min} & &  \delta_1 f_1(\bm{p})  {+} \sum_{\sigma=1}^{\mathcal{E}}\delta_2f_2(\bm{p}_{e\sigma})
 \\
& \: \text{s.t.} & & \eqref{1s}-\eqref{13}, \eqref{21s}-\eqref{23s}
\end{aligned}
\end{pequation}where $\tilde{\bm{u}}=[\tilde{\bm{u}}_1, \ldots, \tilde{\bm{u}}_m]$. Problem \eqref{p2} is coupled through the global power loss objective in \eqref{14} and the global voltage constraint in \eqref{13}. Primal-dual based algorithms \cite{koshal2011multiuser,liu2017decentralized,huo2022two} can solve \eqref{p2} in a decentralized fashion. To this end, we first derive the relaxed Lagrangian function of problem \eqref{p2} where only global constraints are included \vspace{-2mm}
\begin{align}
    \mathcal{L}(\bm{p}, \tilde{\bm{u}}, \bm{\lambda},\bm{\mu}) &= \delta_1 f_1(\bm{p}) {+} \sum_{\sigma=1}^{\mathcal{E}}\delta_2f_2(\bm{p}_{e\sigma}) {+}  \sum_{i=1}^{n}\bm{\lambda}_{i}^{\mathsf{T}} (\boldsymbol{V}_i - \hat{\bm{V}}_{u} )\nonumber\\
    & ~~~ - \sum_{i=1}^{n}\bm{\mu}_{i}^{\mathsf{T}} (\boldsymbol{V}_i - \hat{\bm{V}}_{l} )
\label{Lagrangian}
\end{align}
where $\tilde{\bm{u}} {=} [\tilde{\bm{u}}_1,\ldots,\tilde{\bm{u}}_m]$, $\bm{\lambda} {=} [\bm{\lambda}_1, \ldots, \bm{\lambda}_n]$, $\bm{\mu} {=} [\bm{\mu}_1, \ldots, \bm{\mu}_n]$,  $\bm{\lambda}_{i}$ and $ \bm{\mu}_{i}$ denote the dual variables associated with node $i$. Then, we take the shrunken-primal-dual-subgradient (SPDS) algorithm in \cite{liu2017decentralized} for example to update the  decision variables $\bm{p}_{v\epsilon}$, $\bm{p}_{sk}$, 
$\bm{p}_{e\sigma}$, and $\tilde{\bm{u}}_{j}$ in \eqref{p2}. Let $\mathcal{X}_\imath$ and 
$\mathcal{Y}$ universally denote different primal and dual variables. The SPDS updates the primal variable (decision variable) and dual variable (Lagrange multiplier) by  
\begin{subequations} \label{20}
\begin{align} \mathcal{X}_{\imath}^{(\ell+1)} &{=}\Pi_{\mathbb{X}_{\imath}}\left(\frac{1}{\tau_{\mathcal{X}}} \Pi_{\mathbb{X}_{\imath}}\left(\tau_{\mathcal{X}} \mathcal{X}_{\imath}^{(\ell)}{-}\alpha_{(\imath{,}\ell)} \nabla_{\mathcal{X}_{\imath}} \mathcal{L}\left(\cdot\right)\right)\right) \label{20a}\\ 
\mathcal{Y}^{(\ell+1)} &{=}\Pi_{\mathbb{D}}\left(\frac{1}{\tau_{\mathcal{Y}}} \Pi_{\mathbb{D}}\left(\tau_{\mathcal{Y}} \mathcal{Y}^{(\ell)}{+}\beta_{\ell} \nabla_{\mathcal{Y}} \mathcal{L}\left(\cdot\right)\right)\right) \label{20b}
\end{align}
\end{subequations} $\forall \imath = 1, \ldots,\mathcal{N}$ where $\mathcal{N}$ denotes the number of decision variables, $0<\tau_{\mathcal{X}}, \tau_{\mathcal{Y}} <1$ denote the shrunken parameters for the primal and dual updates, respectively. $\Pi$ denotes the Euclidean projection, $\alpha_{(\imath, \ell)}$ and $\beta_{\ell}$ denote the primal and dual step sizes at the $\ell$th iteration. $\mathbb{X}_{\imath}$  and $\mathbb{D}$ are the feasible sets for the primal variable $\mathcal{X}_{\imath}$ and dual variable $\mathcal{Y}$, respectively. $\nabla \mathcal{L}\left(\cdot\right)$ denotes the first-order subgradient \emph{w.r.t.} the relaxed Lagrangian function. Note that the initial values of the primal or dual variable will not affect the convergence of SPDS.

Then, the subgradients of \eqref{Lagrangian} \emph{w.r.t.} the primal variables $\bm{p}_{v\epsilon}$, $\bm{p}_{sk}$, $\bm{p}_{e\sigma}$,  and $\tilde{\bm{u}}_{j}$ can be calculated, respectively. Due to the page limit,  we only give the calculation of subgradients related to the $\sigma$th ESS as an example. Suppose the $\sigma$th ESS is connected at node $i$, we have 
\begin{align}
\label{19}
\nabla_{\bm{p}_{e\sigma}} \mathcal{L}(\cdot) &= 2 \delta_2 \bm{B}^{\mathsf{T}}\bm{B}\bm{p}_{e\sigma}
    - \frac{2\delta_1}{|V_0|^2}\sum_{l=1}^{n} r_l (\nabla_{\bm{p}_{e\sigma}}\bm{\mathcal{P}}_l) \bm{\mathcal{P}}_l \nonumber\\
    &~~~ + \sum_{\iota=1}^{n} 2\bm{R}_{\iota i} (\bm{\mu}_{\iota} - \bm{\lambda}_{\iota} )
\end{align}
where $r_l$ and $\bm{\mathcal{P}}_l$ denote the resistance and active power flow of the $l$th line, respectively. $\nabla_{\bm{p}_{e\sigma}}\bm{\mathcal{P}}_\iota$ can be calculated based on the distribution network topology and active power injections.  For example, for the 13-node distribution network shown in Fig. \ref{buildings_distribution_network}, the $\sigma$th ESS connected at node 6 has $\nabla_{\bm{p}_{e\sigma}}\bm{\mathcal{P}}_\iota = \bm{\psi} \bm{p}$ where $\bm{p} = [\bm{p}_1^\mathsf{T}, \ldots, \bm{p}_n^\mathsf{T}]^\mathsf{T}$, $\bm{\psi} = [\bm{I},\ldots,\bm{I},2\bm{I},\ldots,2\bm{I}]$, and $\bm{I} \in \mathbb{R}^{T\times T}$ is an identity matrix. Note that $\bm{\psi}$ relates to the upstream lines $\bm{\mathcal{P}}_1$ and $\bm{\mathcal{P}}_6$ that are defined through the network topology and line flow direction. Follow this procedure,  $\nabla_{\bm{p}_{sk}} \mathcal{L}(\cdot)$, $\nabla_{\bm{p}_{v\epsilon}} \mathcal{L}(\cdot)$, and $\nabla_{\tilde{\bm{u}}_{j}} \mathcal{L}(\cdot)$ can be readily derived. 

In what follows, we rewrite the local constraints of each decision variable into a more compact form as 
\begin{subequations}
\begin{align}
\mathbb{P}_\epsilon^v &\triangleq \{ \bm{p}_{v \epsilon} | \: \eqref{1s}, \eqref{2s} \},
\mathbb{P}_k^s \triangleq \{ \bm{p}_{sk} | \: \eqref{pv_limit}\} \nonumber\\
\mathbb{P}_\sigma^e &\triangleq
\{ \bm{p}_{e\sigma} | \: \eqref{charge_limit}, \eqref{state_limit} \}, \mathbb{U}_j \triangleq
\{ \tilde{\bm{u}}_{j} | \: \eqref{21s}, \eqref{22s}, \eqref{17s} \}. \nonumber
\end{align}
\end{subequations}
The feasible set of dual variables is
\begin{align}
\mathbb{D} \triangleq \{ \bm{\lambda}_{i}, \bm{\mu}_{i} | \: \bm{\lambda}_{i}, \bm{\mu}_{i} \geq \bm{0}, \forall i \in \mathbb{N}\}.
\end{align}
Finally, the primal-dual variables can be updated using the derived subgradients and feasible sets by following \eqref{20}. The convergence error of the $j$th aggregator can be calculated by
\begin{align}
 \epsilon_j^{(\ell)} &= \sum_{\epsilon=1}^{\mathcal{V}_j}\| \bm{p}_{v \epsilon}^{(\ell+1)} - \bm{p}_{v \epsilon}^{(\ell)} \|_2 + \sum_{k=1}^{\mathcal{S}_j} \|\bm{p}_{sk}^{(\ell+1)} - \bm{p}_{sk}^{(\ell)}\|_2 \nonumber\\
    & ~~~ + \sum_{\sigma=1}^{\mathcal{E}_j}\|\bm{p}_{e\sigma}^{(\ell+1)} - \bm{p}_{e\sigma}^{(\ell)}\|_2 + \|\tilde{\bm{u}}_{j}^{(\ell+1)} - \tilde{\bm{u}}_{j}^{(\ell)}\|_2
\end{align} where $\mathcal{V}_j$, $\mathcal{S}_j$, and $\mathcal{E}_j$ denote the number of EVs, PVs, and ESSs at building $j$, respectively. Note that the iterative updates in \eqref{20} require the   primal and dual variables to be exchanged in each iteration. In reality, such a setting is usually expensive and challenged by heterogeneous temporal scales of the agents, aggregators, and SO. To resolve this issue, we next introduce the asynchronous update between aggregators and the SO, which can ensure enhanced compatibility and reduce the  computing burden imposed on individual agents. 


\subsection{Asynchronous Updates} \label{Asynchronous_update}

To future improve the communicating compatibility and system resilience  in terms of heterogeneous temporal scales, asynchronous primal update using inexact dual solutions is investigated. The detailed updating procedure is presented via \textbf{Algorithm \ref{alg:1}}.

 \RestyleAlgo{ruled} 
\begin{algorithm}
\caption{Two-level HDC control strategy with asynchronous primal updates} \label{alg:1}
 \algorithmfootnote{$^\star$\hspace{-2mm}$\mod$denotes the modular operation.} 
\textbf{Initialization:} Aggregators initialize decision variables, tolerance $\epsilon_0$,  iteration counter $\ell=0$, and maximum iteration $\ell_{max}$;

\textbf{RC model training:} The $j$th aggregator trains the RC values of the rooms inside the $j$th building using particle swarm optimization (PSO)\cite{cui2019hybrid,cui2022model}, then calculates coefficients $a_j$, $b_j$, and $g_j$ in \eqref{21s};

 \While{$\epsilon_j^{(\ell)} > \epsilon_0$ and $\ell < \ell_{max}$}{
  The $j$th aggregator determines $k_j^{(\ell)}, 1 \leq k_j^{(\ell)} \leq \mathcal{K}$ based on its real-time local computing ability\;
 
 \If{$\ell \mod$\hspace{-2mm} $^\star$  $k_j^{(\ell)} = 0$}{
  The $j$th aggregator updates the decision variables $\bm{p}_{v\epsilon}^{(\ell)} \rightarrow \bm{p}_{v\epsilon}^{(\ell+1)}$, $\bm{p}_{sk}^{(\ell)} \rightarrow \bm{p}_{sk}^{(\ell+1)}$, $\bm{p}_{e\sigma}^{(\ell)} \rightarrow \bm{p}_{e\sigma}^{(\ell+1)}$, $\tilde{\bm{u}}_{j}^{(\ell)} {\rightarrow} \tilde{\bm{u}}_{j}^{(\ell+1)}$ of building $j$;
   
   The $j$th aggregator determines the \textsc{On-Off} states of HVACs using $\tilde{\bm{u}}_{j}^{(\ell+1)}$;

   The $j$th aggregator sends the updated decision variables, i.e., $\bm{p}_{v\epsilon}^{(\ell+1)}$, $\bm{p}_{sk}^{(\ell+1)}$, $\bm{p}_{e\sigma}^{(\ell+1)}$, and $\tilde{\bm{u}}_{j}^{(\ell+1)}$ to the SO\;
   }
   { The SO updates the dual variables $\bm{\lambda}_{i}^{(\ell)} \rightarrow \bm{\lambda}_{i}^{(\ell+1)}$, $ \bm{\mu}_{i}^{(\ell)} \rightarrow \bm{\mu}_{i}^{(\ell+1)}$, $\forall i = 1,\ldots n$ using \eqref{20b}, then sends the updated dual variables $\bm{\lambda}_{i}^{(\ell+1)}$, $\bm{\mu}_{i}^{(\ell+1)}$ $\forall i = 1,\ldots n$ to the aggregators\;
  }
  {$\ell \rightarrow \ell+1$};
 }
\end{algorithm}


 
\textbf{Theorem 1:} (\textit{Inexact dual solutions}) For each update in the dual space, let $k_j, j =1,\ldots,m$, $1\leq  k_j \leq \mathcal{K}$ denote the number of fixed primal updates for the $j$th aggregator. Then, \textbf{Algorithm \ref{alg:1}} converges to the optimal solutions and tolerates maximum $\mathcal{K}\Delta T$ time discrepancy for the updates of the aggregators in the primal space. \hfill $\blacksquare$

In real-world industrial applications, it can be demanding on ensuring synchronous temporal scales across diverse agents, aggregators, and the SO. \textbf{Theorem 1} allows each aggregator to choose its optimal updating time based on its real-time computing ability and communication capacity. Though asynchronous updates can admittedly increase the total number of iterations compared with the case when $\mathcal{K}=1$, the convergence of decision variables and the satisfaction of both global and local constraint are guaranteed. Please refer to \textsc{Appendix} \ref{Theorem1_proof} for the proof of \textbf{Theorem 1}.



\section{Simulation Results}

\begin{figure*}%
\vspace{-4mm}
\centering
\begin{subfigure}[t]{.47\columnwidth} \centering
\includegraphics[width=1.07\columnwidth,trim = 0mm 3mm 0mm 0mm, clip]{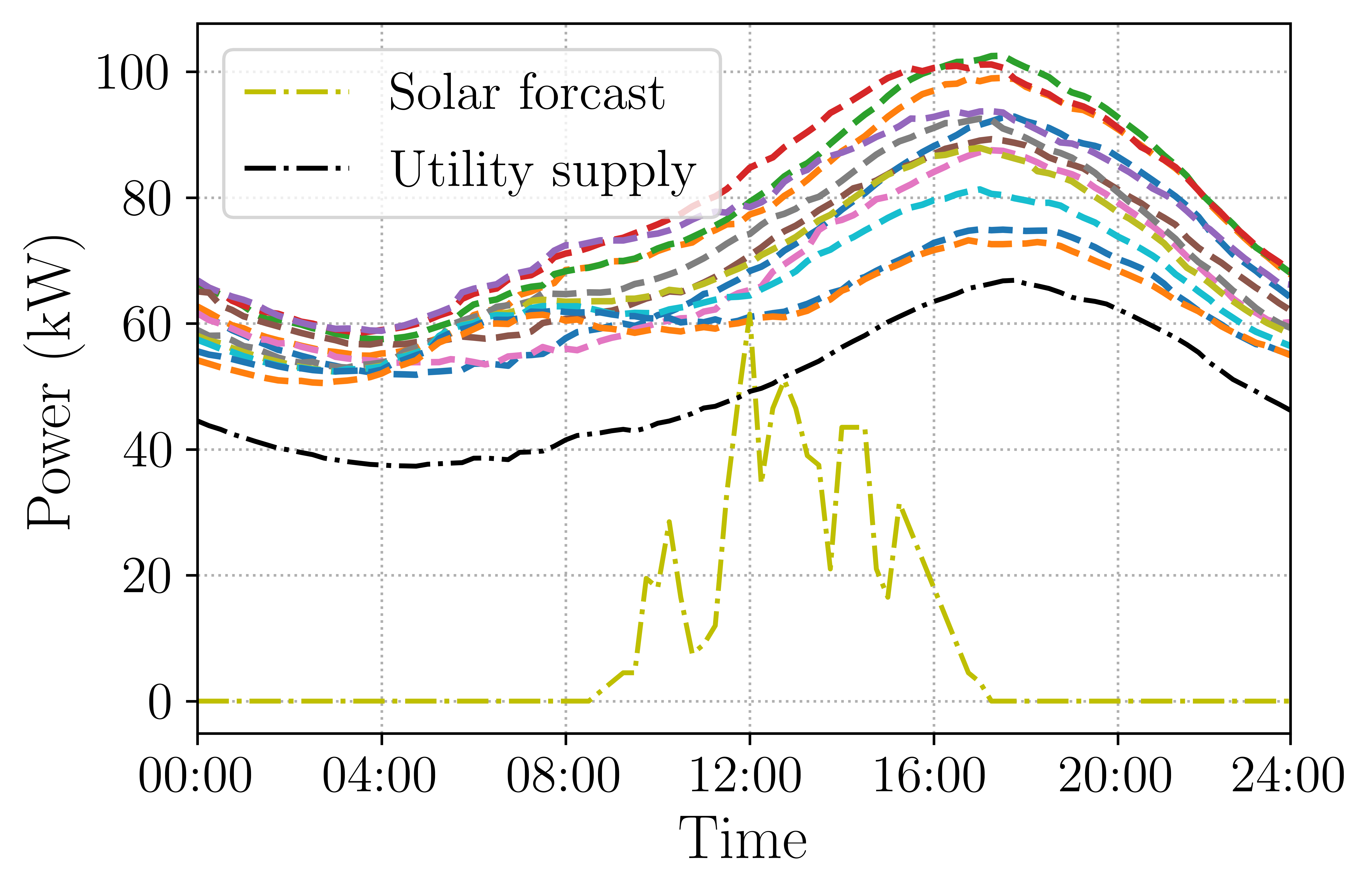}%
\caption{Baseline loads, solar prediction, and the utility supply.}%
\label{subfig3a}%
\end{subfigure}\hfill
\begin{subfigure}[t]{.49\columnwidth} \centering
\includegraphics[width=1.03\columnwidth, trim = 0mm 3mm 0mm 0mm, clip]{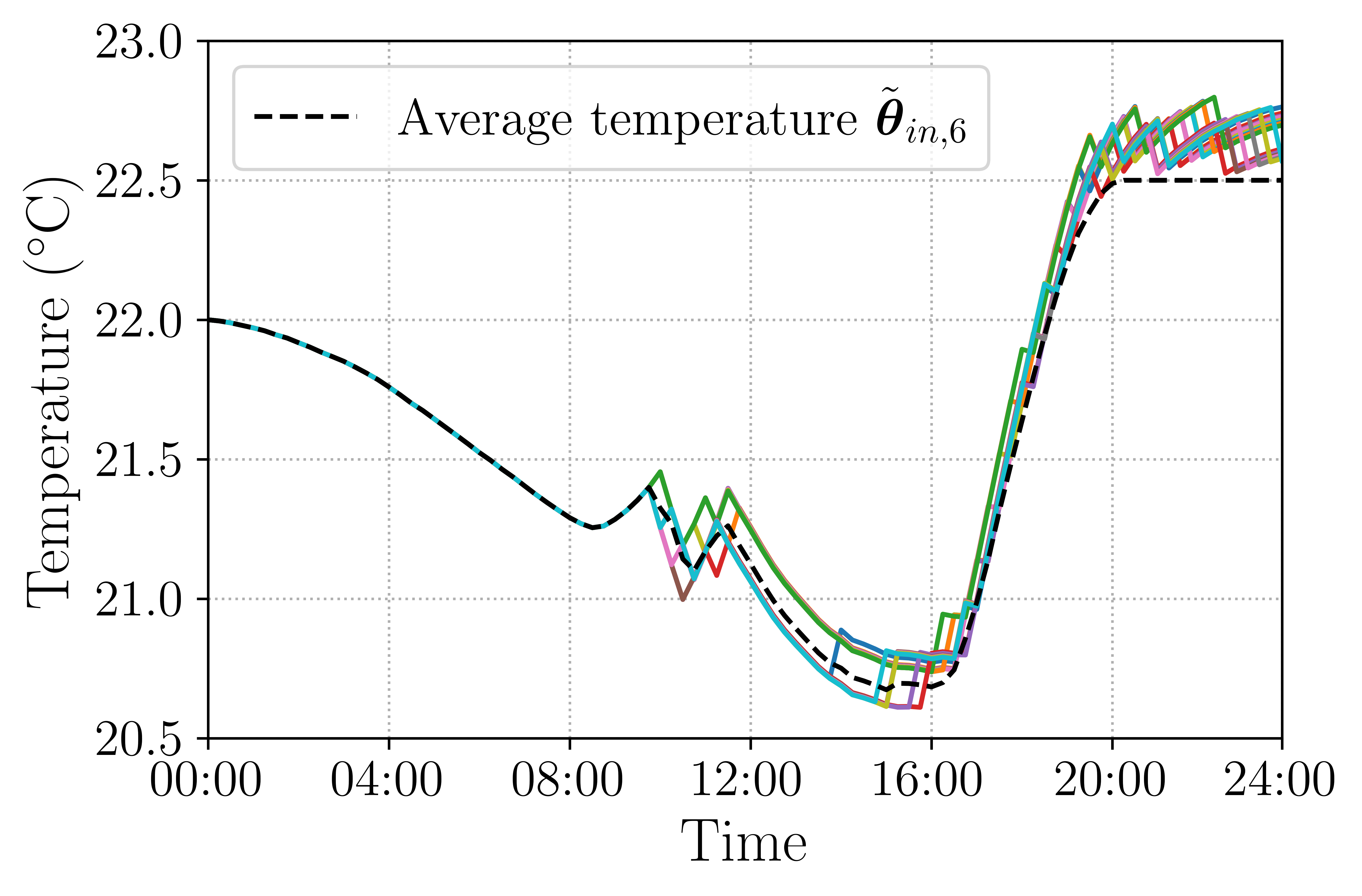}%
\caption{Individual room temperatures of Building 6 (Solid lines).}%
\label{subfig3b}%
\end{subfigure}\hfill%
\begin{subfigure}[t]{.47\columnwidth}
\centering
\includegraphics[width=1.02\columnwidth, trim = 0mm 3mm 0mm 0mm, clip]{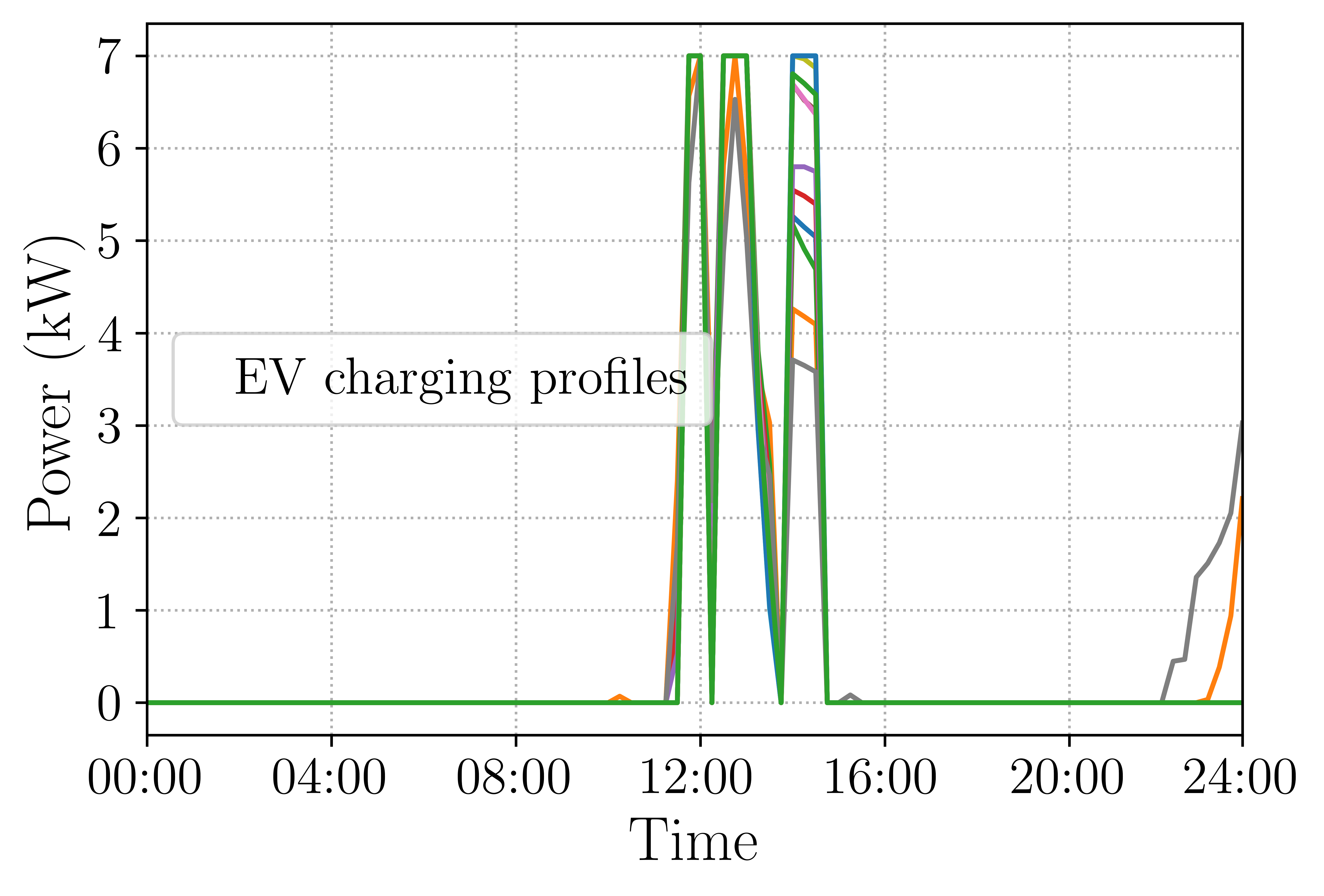}%
\caption{EVs' charging schedules.}%
\label{subfig3c}%
\end{subfigure} \hfill%
\begin{subfigure}[t]{.47\columnwidth}
\centering
\includegraphics[width=1.05\columnwidth, trim = 0mm 3mm 0mm 0mm, clip]{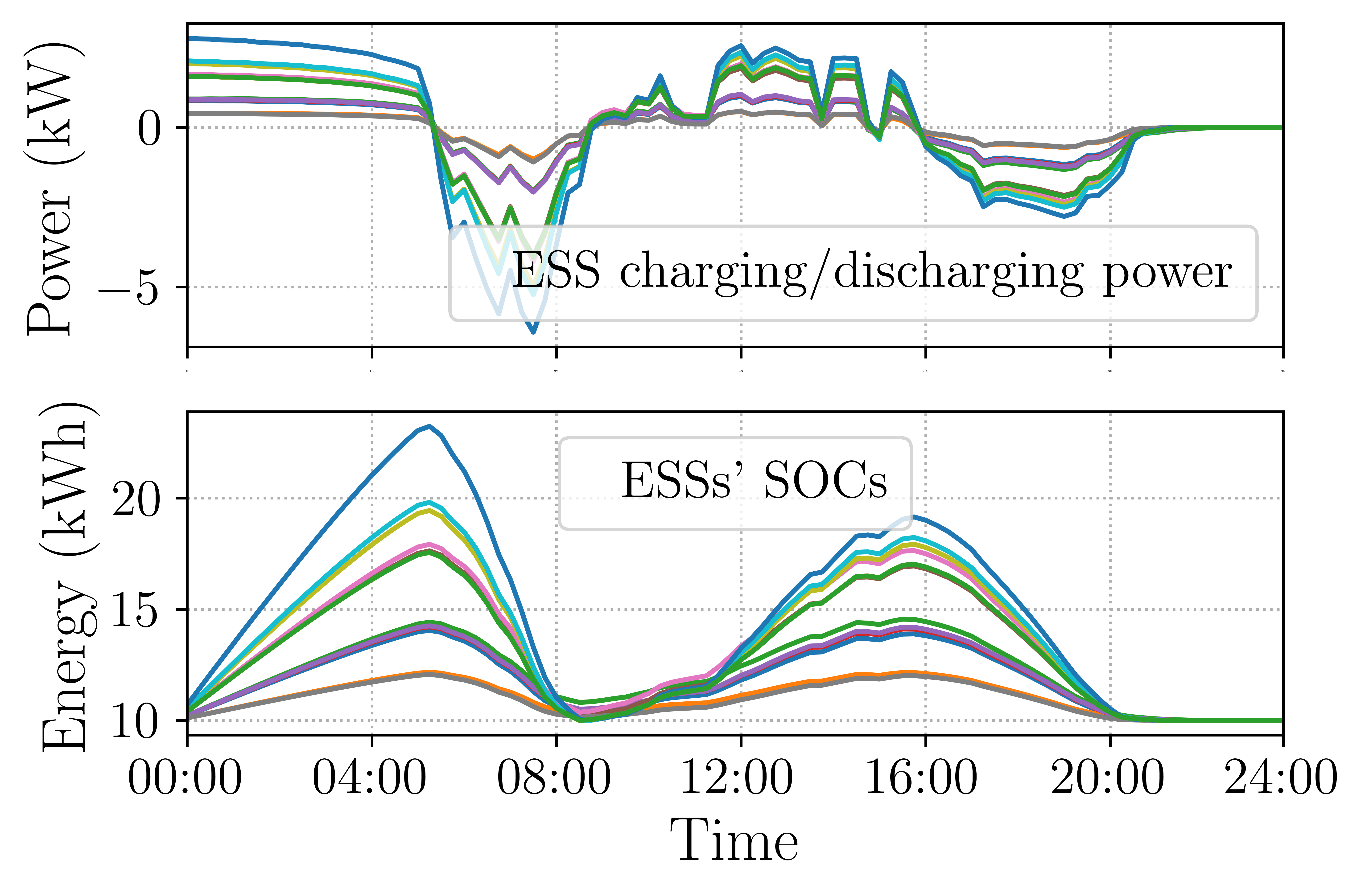}%
\caption{ESSs' charging/discharging behaviors and SOCs.}%
\label{subfig3d}%
\end{subfigure}%

\caption{The optimized DER and HVAC schedules using the proposed HDC algorithm.}
\label{figabcd}
\vspace{-2mm}
\end{figure*}


\begin{figure}
    \vspace{-3mm}
    \centering
\includegraphics[width=0.95\columnwidth,trim = 0mm 3mm 70mm 180mm, clip]{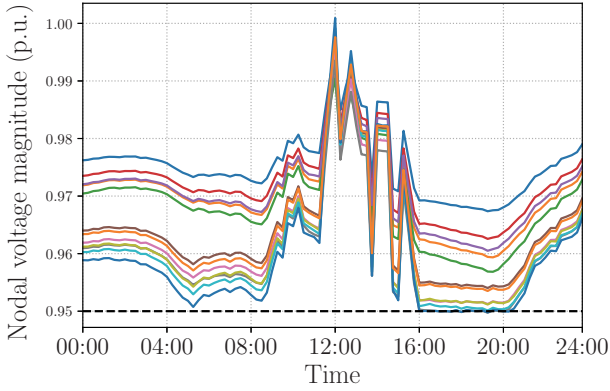}    
\caption{Nodal voltage magnitudes.}
    \label{fig_4}
    \vspace{-5mm}
\end{figure}

Historic data from Building 3147 were used to train and validate the RC model. The RC parameters of all the rooms were trained by using the PSO routine, more details regarding the training process can be found in  \textsc{Appendix} \ref{RC_PSO}. Without loss of generality, we assume the buildings in the distribution network share the same weather conditions, i.e., solar radiation and outdoor temperature.

To demonstrate the grid-level potentials of the proposed HDC algorithm, we assume each node in the IEEE 13-node test feeder (Fig. \ref{buildings_distribution_network}) is connected with one building resembling Building 3147. In total 8 rooftop PVs are assumed to be connected at each building, resulting in around 50 kW peak solar power generation, and the PV data is collected from the DECC lab on a sunny day on September 6th 2021. Besides, one ESS and one EV are assumed to be connected at each building. The ESSs' capacity limits are calibrated to be $\hat{\bm{p}}_e^{cl} = 10$ kWh and $\hat{\bm{p}}_e^{cu} = 60$ kWh, the maximum charging/discharging rates are $\hat{\bm{p}}_{e}^{l} = -15$ kW and $\hat{\bm{p}}_{e}^{u} = 15$ kW, respectively. The EV is charged by a level-2 EV charger with a maximum charging power 7.6 kW and the charging demand of all the EVs randomly distribute in $[12,16]$ kWh. The lower and upper voltage limits are set to be $\hat{\bm{V}}_{l} = 0.95 \bm{V}_0$ and $\hat{\bm{V}}_{u} = 1.05 \bm{V}_0$, respectively. The time interval is chosen to be $\Delta T = 15$ min. The aggregator can control $N_j = 20$ HVACs in each building, and the aggregated building parameters are trained and uniformly set to be $a_j = 0.98$, $b_j = 0.02$, and $g_j = -0.2$, respectively. The HVACs' rated power is uniformly chosen as $P^{r}= 1 ~\text{kW}$. The temperature comfort zone is set to be $[20^{\circ}\text{C}, 22.5^{\circ}\text{C}]$.  To test the efficacy of asynchronous primal updates, let the aggregators of Building 2 and Building 8 update their decision variables every $2 \Delta T$ and $4 \Delta T$, respectively. 

Fig. \ref{figabcd} shows the scheduling and control results of all DERs and HVACs using the proposed HDC algorithm. Specifically, Fig. \ref{subfig3a} presents the baseline loads of in total 12 buildings, the solar prediction, and the utility supply, respectively. Fig. \ref{subfig3b} shows the temperature control of Building 6 where the black dashed line gives the aggregated temperature calculated with the aggregated model. The solid lines give the individual room temperatures whose absolute deviations from the aggregated temperature were less than $0.5^{\circ}\text{C}$. The EVs' charging schedules are shown in Fig. \ref{subfig3c}, where all EVs charged at the fastest rate at noon as a result of sufficient solar radiation. As shown in Fig. \ref{subfig3d}, ESSs stored solar energy and charged mostly during the noon time, then discharged the stored energy starting from $16:00$ in the afternoon when the solar power decreases while the load demand increases as shown in Fig. \ref{subfig3a}, and the state of charges (SOCs) of all ESSs are maintained above the limits through the day. The active power loss $f_1(\bm{p})$ and battery degradation cost $\sum_{\sigma=1}^{\mathcal{E}}f_2(\bm{p}_{e\sigma})$ obtained in the IEEE 13-node system are
0.18 kW, and 229 (kW)$^2$, respectively. Finally, Fig. \ref{fig_4} gives the nodal voltage magnitudes of 12 nodes, only the lower voltage bound is active due to a large amount of loads, and all the voltage magnitudes are above the lower voltage limit which is represented by the dotted black line.

\section{Conclusions}

This paper developed a scalable, efficient, and compatible control framework for GEBs to provide both building-level and grid-level services. The control problem was formulated into a constrained optimization problem consisting of global and local objectives and constraints. We proposed an HDC algorithm to solve the formulated problem in a hybrid decentralized-centralized fashion. The HDC structure can scale to a large number of grid-connected devices and buildings and adapt to asynchronous updates for different entities under heterogeneous temporal scales. We verified the efficiency and efficacy of the proposed HDC algorithm via real-world collected data and conducted simulations based on an office building located at ORNL. Future work in this direction includes investigating real-time building-level scheduling as well as control strategies. 

 \begin{appendices}
\section{Proof of \textbf{Theorem 1}}\label{Theorem1_proof} 

For the SPDS updates in \eqref{20}, consider the execution of $\mathcal{K}$ primal updates for every dual update, we have 
\begin{subequations} \label{20ss}
\begin{align} \mathcal{X}^{(k+1)}[\mathcal{Y}^{(\ell)}] &=\Pi_{\mathbb{X}}\left(\frac{1}{\tau_{\mathcal{X}}} \Pi_{\mathbb{X}}\left(\tau_{\mathcal{X}} \mathcal{X}^{(k)}[\mathcal{Y}^{(\ell)}]\right.\right. \nonumber\\
&~~ -\left.\left.\alpha\nabla_{\mathcal{X}} \mathcal{L}(\mathcal{X}^{(k)}[\mathcal{Y}^{(\ell)}], \mathcal{Y}^{(\ell)})\vphantom{\mathcal{X}^{(k)}}\right)\vphantom{\frac{1}{\tau_{\mathcal{X}}}}\right)\label{20ass}\\ 
\mathcal{Y}^{(\ell+1)} &=\Pi_{\mathbb{D}}\left(\frac{1}{\tau_{\mathcal{Y}}} \Pi_{\mathbb{D}}\left(\tau_{\mathcal{Y}} \mathcal{Y}^{(\ell)}\right.\right.\nonumber\\
&~~+\left.\left.\beta \nabla_{\mathcal{Y}} \mathcal{L}(\mathcal{X}^{(\mathcal{K})}[\mathcal{Y}^{(\ell)}], \mathcal{Y}^{(\ell)})\vphantom{\mathcal{X}^{(k)}}\right)\vphantom{\frac{1}{\tau_{\mathcal{X}}}}\right) \label{20bss}
\end{align}
\end{subequations} $\forall k=1,\ldots,\mathcal{K}$, where $\mathcal{X}^{(k)}[\mathcal{Y}^{(\ell)}]$ denotes the primal solution obtained by SPDS using dual solution $\mathcal{Y}^{(\ell)}$. 
For simplicity, we omit the subscript $\imath$ of $\mathcal{X}_{\imath}$, let $\mathcal{X}^{(k)}[\ell] \triangleq \mathcal{X}^{(k)}[\mathcal{Y}^{(\ell)}]$, and adopt universal primal and dual step sizes $\alpha$ and $\beta$ through the proof. 

Let $\zeta^{*} \triangleq [{\mathcal{X}[\ell]}^{\mathsf{T}}, {\lambda^{(\ell)}}^{\mathsf{T}}]^{\mathsf{T}}$ denote the primal-dual solutions at the $\ell$th iteration and $\zeta^{(k)} \triangleq [{\mathcal{X}^{(k)}[\ell]}^{\mathsf{T}}, {\lambda^{(\ell)}}^{\mathsf{T}}]^{\mathsf{T}}$. By using the nonexpansive property of the projection $\Pi_{\mathbb{X}}(\cdot)$, we have 
\begin{align}
&\left\|\mathcal{X}^{(K)}[\ell]-\mathcal{X}[\ell]\right\|_{2}^{2} \nonumber \\
& ~\leq\frac{1}{\tau_{\mathcal{X}}^{2}}\left\|\mathcal{X}^{(K-1)}[\ell]-\mathcal{X}[\ell]\right\|_{2}^{2}+\frac{\alpha^{2}}{\tau_{\mathcal{X}}^{2}}\|\Phi_{1}(\zeta^{(K-1)})-\Phi_{1}\left(\zeta^{*})\right\|_{2}^{2} \nonumber\\
& \quad -2 \frac{\alpha}{\tau_{\mathcal{X}}^{2}}(\Phi_{1}(\zeta^{(K{-}1)}){-}\Phi_{1}(\zeta^{*}))^{\mathsf{T}}(\mathcal{X}^{(K{-}1)}[\ell]{-}\mathcal{X}[\ell])
\label{26sp}
\end{align} where the mapping $\Phi_{1}(\mathcal{X}, \mathcal{Y}) = \nabla_{\mathcal{X}}\mathcal{L}(\mathcal{X}, \mathcal{Y})+\frac{1-\tau_{\mathcal{X}}}{\alpha} \mathcal{X}$. 

We first deal with the second term on the right-hand side of \eqref{26sp}, we have 
\begin{align} \label{23s}
& \left\Vert \Phi_{1}(\zeta^{(K-1)})-\Phi_{1}(\zeta^{*})\right\rVert_{2}\nonumber\\
& \quad \leq  \Big\Vert \nabla_{\mathcal{X}} \mathcal{L}(\mathcal{X}^{(K-1)}[\ell], \mathcal{Y}^{(\ell)}) - \nabla_{\mathcal{X}} \mathcal{L}(\mathcal{X}[\ell], \mathcal{Y}^{(\ell)})\Big\Vert_2 \nonumber\\
  & \qquad +\frac{1-\tau_\mathcal{X}}{\alpha}\Big\Vert \mathcal{X}^{(K-1)}[\ell]-\mathcal{X}[\ell]\Big\Vert_2 \nonumber\\
  & \quad \leq \left\Vert \nabla_\mathcal{X} G(\mathcal{X}^{(K-1)}[\ell]){-} \nabla_\mathcal{X} G(\mathcal{X}[\ell]) \right\Vert_2 {+} \rho\left\Vert \mathcal{X}^{(K-1)}[\ell]{-}\mathcal{X}[\ell]\right\rVert_2 \nonumber\\
  & \qquad + \left\Vert \nabla_\mathcal{X}d^{\mathsf{T}}(\mathcal{X}^{(K-1)}[\ell])\mathcal{Y}^{(\ell)}- \nabla_\mathcal{X}d^{\mathsf{T}}(\mathcal{X}[\ell])\mathcal{Y}^{(\ell)}\right\Vert_2 \nonumber \\
 & \qquad +\frac{1-\tau_\mathcal{X}}{\alpha}\Big\Vert \mathcal{X}^{(K-1)}[\ell]-\mathcal{X}[\ell]\Big\Vert_2 \nonumber\\
& \quad \leq L_{\nabla G} \Big\Vert \mathcal{X}^{(K-1)}[\ell]-\mathcal{X}[\ell]\Big\Vert_2 + \rho\left\Vert \mathcal{X}^{(K-1)}[\ell]{-}\mathcal{X}[\ell]\right\rVert_2 \nonumber\\
& \qquad +\frac{1-\tau_\mathcal{X}}{\alpha}\Big\Vert \mathcal{X}^{(K-1)}[\ell]-\mathcal{X}[\ell]\Big\Vert_2 \nonumber\\
& \quad = L_{\Phi_1} \Big\Vert \mathcal{X}^{(K-1)}[\ell]-\mathcal{X}[\ell]\Big\Vert_2 
\end{align}
where $L_{\nabla G}$ denotes the Lipschitz constant
 of $\nabla_{\mathcal{X}} G(\mathcal{X})$, $L_{\Phi_1} \triangleq L_{\nabla G} + \rho + \frac{1-\tau_\mathcal{X}}{\alpha}$, and \eqref{23s} comes from the fact that $\nabla_\mathcal{X}d^{\mathsf{T}}(\mathcal{X}^{(K-1)}[\ell])\mathcal{Y}^{(\ell)}- \nabla_\mathcal{X}d^{\mathsf{T}}(\mathcal{X}[\ell])\mathcal{Y}^{(\ell)}=0$.

Therefore, we can readily obtain 
\begin{equation}\label{28sp}
 \left\Vert \Phi_{1}(\zeta^{(K-1)}){-}\Phi_{1}(\zeta^{*})\right\rVert_{2}^2 {\leq} L_{\Phi_1}^2 \Big\Vert \mathcal{X}^{(K-1)}[\ell]-\mathcal{X}[\ell]\Big\Vert_2^2.
\end{equation}

Then, evaluating the second term on the right-hand side of \eqref{26sp}, we have 
\begin{align} \label{29sp}
& \left( \Phi_{1}(\zeta^{(K-1)})-\Phi_{1}(\zeta^{*})\right)^{\mathsf{T}}  (\mathcal{X}^{(K{-}1)}[\ell]{-}\mathcal{X}[\ell])\nonumber\\
& \quad =  \left( \nabla_{\mathcal{X}} \mathcal{L}(\mathcal{X}^{(K{-}1)}[\ell], \mathcal{Y}^{(\ell)}) {-} \nabla_{\mathcal{X}} \mathcal{L}(\mathcal{X}[\ell], \mathcal{Y}^{(\ell)})\right)^{\mathsf{T}} (\mathcal{X}^{(K{-}1)}\nonumber\\
& \qquad -\mathcal{X}[\ell]) +\frac{1-\tau_\mathcal{X}}{\alpha}\Big\Vert \mathcal{X}^{(K-1)}[\ell]-\mathcal{X}[\ell]\Big\Vert_2^2  \nonumber\\
& \quad = \left( \nabla_\mathcal{X} G(\mathcal{X}^{(K-1)}[\ell]){-} \nabla_\mathcal{X} G(\mathcal{X}[\ell]) \right)^{\mathsf{T}}  (\mathcal{X}^{(K{-}1)}[\ell]{-}\mathcal{X}[\ell]) \nonumber\\
&\qquad + \rho\left\Vert \mathcal{X}^{(K-1)}[\ell]{-}\mathcal{X}[\ell]\right\Vert_2^2 \nonumber\\
& \qquad +\left(\nabla_\mathcal{X}d^{\mathsf{T}}(\mathcal{X}^{(K{-}1)}[\ell])\mathcal{Y}^{(\ell)}{-}\nabla_\mathcal{X}d^{\mathsf{T}}(\mathcal{X}[\ell])\mathcal{Y}^{(\ell)}\right)^{\mathsf{T}} \hspace{-2mm} (\mathcal{X}^{(K{-}1)}\nonumber\\
& \qquad -\mathcal{X}[\ell]) +\frac{1-\tau_\mathcal{X}}{\alpha}\Big\Vert \mathcal{X}^{(K-1)}[\ell]-\mathcal{X}[\ell]\Big\Vert_2^2.
\end{align}
Since the function $G(\mathcal{X})$ is convex, we have $\left( \nabla_\mathcal{X} G(\mathcal{X}^{(K-1)}[\ell]){-} \nabla_\mathcal{X} G(\mathcal{X}[\ell]) \right)^{\mathsf{T}}  (\mathcal{X}^{(K{-}1)}[\ell]{-}\mathcal{X}[\ell]) \geq 0$. Therefore, \eqref{29sp} becomes 
\begin{align} \label{30sp}
& \left( \Phi_{1}(\zeta^{(K-1)})-\Phi_{1}(\zeta^{*})\right)^{\mathsf{T}}  (\mathcal{X}^{(K{-}1)}[\ell]{-}\mathcal{X}[\ell])\nonumber\\
& \quad \geq (\rho + \frac{1-\tau_\mathcal{X}}{\alpha}) \left\Vert \mathcal{X}^{(K-1)}[\ell]{-}\mathcal{X}[\ell]\right\Vert_2^2.
\end{align}

Substituting \eqref{28sp} and \eqref{30sp} into \eqref{26sp}, we have 
\begin{equation}
\left\|\mathcal{X}^{(K)}[\ell]-\mathcal{X}[\ell]\right\|_{2}^{2} \leq \varsigma\left\|\mathcal{X}^{(K-1)}[\ell]-\mathcal{X}[\ell]\right\|_{2}^{2}
\label{31sp}
\end{equation}
where $\varsigma = \frac{1}{\tau_{\mathcal{X}}^{2}} {+} \frac{\alpha^{2}}{\tau_{\mathcal{X}}^{2}}L_{\Phi_1}^2 -2 \frac{\alpha}{\tau_{\mathcal{X}}^{2}}(\rho + \frac{1-\tau_\mathcal{X}}{\alpha})$. Recursively using \eqref{31sp}, we have 
\begin{equation}
\left\|\mathcal{X}^{(K)}[\ell]-\mathcal{X}[\ell]\right\|_{2}^{2} \leq \varsigma^{K-1}\left\|\mathcal{X}^{(1)}[\ell]-\mathcal{X}[\ell]\right\|_{2}^{2}
\label{32sp}
\end{equation}
Finally, by selecting appropriate $\alpha$, $\tau_{\mathcal{X}}$, and $\rho$, $0<\varsigma <1$ can be satisfied. Then, we have $\mathcal{X}^{(K)}[\ell]  \rightarrow \mathcal{X}[\ell]$ as $K\rightarrow \infty$, therefore the convergence of the sequence $\mathcal{X}^{(k)}[\ell], k=1,\ldots,\mathcal{K}$ is proved. Note that the primal solution $\mathcal{X}^{(K)}[\ell]$ will converge to $\mathcal{X}[\ell]$ \emph{w.r.t.} each dual solution $\mathcal{Y}^{\ell}$ by resembling the convergence of SPDS and a sufficient small $\varsigma$ can be tuned to guarantee the convergence of \eqref{20ss} with $\mathcal{K}$ given.

\section{RC Parameters Trained by PSO}\label{RC_PSO} 

The training process aims at finding the optimal values of the undetermined RC parameters, including the thermal capacitances of the exterior walls, indoor air, and internal mass; the thermal resistances of exterior walls, window, and internal mass, and the convection fractions, i.e., $\mathbb{P}\triangleq \{C_{w},C_{in}, C_{m}, R_{w1}, R_{w2}, R_{win}, R_m, Sp_{1}, Sp_{2}, Sp_{3} \}$. To find the optimal set of RC parameters, we formulate the RC parameter searching process as a nonlinear optimization problem with the following objective function
\begin{equation} 
\mathcal{J}(\mathbb{P})  =\sqrt{\frac{\sum_{t=1}^{T}\left(\bar{\theta}_{in}(t) - \hat{\theta}_{in}(t)\right)^{2}}{T-1}}
\label{R10}
\end{equation}
where $\bar{\theta}_{in}(t)$ and $\hat{\theta}_{in}(t)$ denotes the mean indoor air temperature from EnergyPlus\texttrademark $~$ and the mean indoor air temperature predicted by the developed RC model at time $t$, respectively.

Then, we identified the parameters using PSO. A population of candidate solutions was generated first, then 
moved around in the search space. The solutions were updated and optimized iteratively until a given measure of quality is satisfied. The objective function defined in \eqref{R10} was used as the accuracy indicator to evaluate the performance of the developed RC model. In specific, we set the search ranges for $C_{w}$, $C_{in}$, $R_{w1}$, $R_{w2}$, $R_{win}$, $R_m$ according to the recommended parameters values of buildings in Zone 4 from IECC \cite{IECC2006}, $C_{m}$ is assumed to be in the range of 100–450 (kJ/Km$^2$) \cite{dominguez2010uncertainty}. The lower and upper limits for each R and C are 1/3 and 3 times of the estimated values respectively. The ranges of $Sp_1$, $Sp_2$ and $Sp_3$ are estimated based on experience.

\end{appendices}

\bibliographystyle{IEEEtran}

\bibliography{bibliography}

\end{document}